\documentclass[preprint,12pt]{elsarticle}

\usepackage{amssymb}
 \usepackage{graphicx}

\usepackage[cp1251]{inputenc}
\usepackage{epsfig}

\usepackage{amsthm}

\newcommand{\sysn}{\left\{\begin{array}{rcl}}
\newcommand{\sysk}{\end{array}\right.}

\newtheorem{theorem}{Theorem}[section]

\theoremstyle{definition}

\newtheorem{corollary}[theorem]{Corollary}

\journal{...}

\begin{document}

\title{Every metric space of weight $\lambda=\lambda^{\aleph_0}$ admits a condensation onto
a Banach space}

\author[affil1]{Alexander V. Osipov}
\address[affil1]{Krasovskii Institute of Mathematics and Mechanics, \\ Ural Federal
 University, Ural State University of Economics, Yekaterinburg, Russia}
\ead[affil1]{OAB@list.ru}

\author[affil2]{Evgenii G. Pytkeev}
\address[affil2]{Krasovskii Institute of Mathematics and Mechanics, Yekaterinburg, Russia}
\ead[affil2]{pyt@imm.uran.ru}

\begin{abstract} In this paper, we have proved that for each
cardinal number $\lambda$ such that $\lambda=\lambda^{\aleph_0}$ a
metric space of weight $\lambda$ admits a bijective continuous
mapping onto a Banach space of weight $\lambda$. Then, we get that
every metric space of weight continuum admits a bijective
continuous mapping onto the Hilbert cube. This resolves the famous
Banach's Problem (\,{\it when does a metric (possibly Banach)
space $X$ admit a bijective continuous mapping onto a compact
metric space ?}) in the class of metric spaces of weight
continuum.

Also we get that every metric space of weight
$\lambda=\lambda^{\aleph_0}$ admits a bijective continuous mapping
onto a Hausdorff compact space. This resolves the Alexandroff
Problem ({\it when does a Hausdorff space $X$ admit a bijective
continuous mapping onto a Hausdorff compact space?}) in the class
of metric spaces of weight $\lambda=\lambda^{\aleph_0}$.

\end{abstract}

%\tnotetext[label1]{The research has been supported by .}


\maketitle %%

% \linenumbers


\section{Intoduction}

In 1935, Stefan Banach wrote down Problem 1 in the Scottish Book:
{\it when does a metric (possibly Banach) space $X$ admit a
condensation (i.e. a bijective continuous mapping) onto a
compactum (= compact metric space)}?

The answer to this problem for Banach spaces can be found in
\cite{bp,banakh,os,pytk}.

Independently to S. Banach the problems concerning condensations
were posed by P.S. Alexandroff (\cite{arh}, Problem 2.7): {\it
when does a Hausdorff space $X$ admit a condensation onto a
Hausdorff compact space?}

\medskip

The main result of this paper is the following theorem that
answers Banach's Problem in the class of metric spaces of weight
continuum and Alexandroff's Problem in the class of metric spaces
of weight $\lambda=\lambda^{\aleph_0}$.

\begin{theorem} Let $\lambda$ be a cardinal number such that
$\lambda=\lambda^{\aleph_0}$. If $(X,\tau)$ is a metric space of
weight $\lambda$ then $\tau$ is the supremum of two topologies
$\tau_1$ and $\tau_2$, where $(X,\tau_i)$ is homeomorphic to a
Banach space of weight $\lambda$ for $i=1,2$.
\end{theorem}

The motivation and idea of proving the main result of this paper
follows from the result of Theorem 1 in \cite{pytk}: {\it Every
separable absolute Borel space $X$ condenses onto the Hilbert
cube, whenever $X$ is not $\sigma$-compact}.

\section{Main definitions and notation}

Recall that a topology $\tau$ on $X$ is the supremum of topologies
$\tau_1$ and $\tau_2$ on $X$ if it is the coarsest topology on $X$
that is finer that $\tau_i$ for every $i=1,2$.

Let $\lambda$ be an infinite cardinal, $S$ a set of cardinality
$\lambda$, and let $I=[0,1]$ be the closed unit interval. Define
an equivalence relation $E$ on $I\times S$ by $(x, \alpha)E(y,
\beta)$ if either $x=0=y$ or $(x,\alpha)=(y,\beta)$. Let
$H(\lambda)$ be the set of all equivalence classes of $E$; in
other words, $H(\lambda)$ is the quotient set obtained from
$I\times S$ by collapsing the subset $\{0\}\times S$ to a point.
For each $x\in I$ and each $\alpha\in S$, $\langle
x,\alpha\rangle$ denotes the element of $H(\lambda)$ corresponding
to $(x, \alpha)\in I\times S$. There is the topology $\tau_{\rho}$
induced from the metric $\rho$ on $H(\lambda)$ defined by
$\rho(\langle x,\alpha\rangle, \langle y,\beta\rangle)=|x - y|$ if
$\alpha=\beta$ and $\rho(\langle x,\alpha\rangle, \langle
y,\beta\rangle)=x+y$ if $\alpha\neq \beta$. The set $H(\lambda)$
with this topology is called {\it the metrizable hedgehog of
spininess $\lambda$} and is often denoted by
$(J(\lambda),\tau_{\rho})$ (\cite{H10}, 4.1.5). The space
$(J(\lambda),\tau_{\rho})$ is a complete, non-compact, metric
space of weight $\lambda$. Moreover, the product
$J(\lambda)^{\aleph_0}$ is a universal space of metrizable spaces
of weight $\lambda$. Every Banach space of weight
$\lambda\geq\aleph_0$ is homeomorphic to $J(\lambda)^{\aleph_0}$
\cite{tur}.

Fix $d\in S$. The second topology $\tau_d$ on $H(\lambda)$ is the
topology generated by the base of the neighborhood system:

$\{O(z)\in \tau_{\rho}\}$ for each $z\in J(\lambda)\setminus
((0,1]\times \{d\})$ and

$\{O(K,\epsilon, z)=(x-\epsilon,x+\epsilon)\times (S\setminus K):
K\in [S]^{<\omega}, d\not\in K, \epsilon>0\}$ for $z=(x,d)\in
(0,1]\times \{d\}$.

The set $H(\lambda)$ with this topology is denoted by
$(J(\lambda),\tau_{d})$. The space $(J(\lambda),\tau_{d})$ is a
Hausdorff compact space. Note that if  $d_1,d_2\in S$ and $d_1\neq
d_2$ then $\tau_{\rho}=\sup \{\tau_{d_1},\tau_{d_2}\}$.

 For the terms and symbols that we do not define follow
\cite{H10}.

\section{Proof of Theorem 1.1}

\begin{proof} Let $(X,\tau)$ be a metric space with the metric $\rho$ of
weight $\lambda$. Because $\lambda=\lambda^{\aleph_0}$, $\lambda$
is not sequential (K\"{o}nig's theorem; see e.g. (\cite{sirp}, p.
181)), so that (\cite{stone}, Th. 2.7) $X$ has a metrically
discrete subset $D$ of cardinality $\lambda$ ($D$ is metrically
discrete if there is $\epsilon>0$ such that every two distinct
points $x,y\in D$ satisfy $\rho(x,y)\geq \epsilon$). Note that any
subset of $D$ is closed subset of $X$.

Let $D=\bigcup\limits_{i=1}^{4}D_i$ where $D_i\cap D_j=\emptyset$
for $i\neq j$ and $|D_i|=\lambda$ for $i,j\in \{1,2,3,4\}$.

Let $V_i=\{x\in X: \rho(x,D_i)<\frac{\epsilon}{3}\}$. Clearly
$\overline{V_i}\cap \overline{V_j}=\emptyset$ for $i\neq j$ and
$D_i\subseteq V_i$ for $i\in \{1,2,3,4\}$.

Let $S$ be a set of cardinality $\lambda$ and fix distinct points
$d_0, d_1, d_2\in S$. Then $I'=((0,1]\times \{d_0\})\cup
\{0\}\subset J(\lambda)$ is homeomorphic to $I=[0,1]$.

We will construct two condensations $f,g: X\rightarrow
(J(\lambda)\times J(\lambda)^{\aleph_0}\times
J(\lambda)^{\aleph_0})\cong J(\lambda)^{\aleph_0}$ where
$J(\lambda)$ is the metrizable hedgehog of spininess $\lambda$
(i.e. $J(\lambda)=(J(\lambda),\tau_{\rho})$) such that no matter
which topology $\tau_{\rho}$, $\tau_{d_1}$ or $\tau_{d_2}$ is
considered in the space $J(\lambda)$, the properties are
satisfied:

I) $f\upharpoonright (X\setminus (V_3\cup V_4))$ is a
homeomorphism,

II) $g\upharpoonright (X\setminus (V_1\cup V_2))$ is a
homeomorphism,

III) $\overline{f(V_4)}\cap \overline{f(V_1)}=\emptyset$ and
$\overline{f(V_3)}\cap \overline{f(V_2)}=\emptyset$,

IV) $\overline{g(V_2)}\cap \overline{g(V_4)}=\emptyset$ and
$\overline{g(V_1)}\cap \overline{g(V_3)}=\emptyset$.

\medskip

{\it Construction of $f$}.

Since $X$ is a normal space there is a continuous mapping
$\varphi_0: X\rightarrow I'$ such that
$\varphi_0(\overline{V_1}\cup \overline{V_3})=0$ and
$\varphi_0(\overline{V_2}\cup \overline{V_4})=1$. Let
$F_1=\varphi_0^{-1}([0,\frac{1}{2}])$ and
$F_2=\varphi_0^{-1}([\frac{1}{2},1])$.

By Kowalsky's Theorem (\cite{kow} or (\cite{H10}, Th. 4.4.9)), the
space $X\setminus(V_3\cup V_4)$ is embedded in
$(J(\lambda)\setminus ((0,1]\times \{d_1\cup d_2\}))^{\aleph_0}$.
Let $\varphi_1: X\setminus(V_3\cup V_4)\rightarrow
(J(\lambda)\setminus ((0,1]\times \{d_1\cup d_2\}))^{\aleph_0}$ be
a embedding. Then $\varphi_2: X\setminus(V_3\cup V_4)\rightarrow
I'\times J(\lambda)^{\aleph_0}\times \{a^0\}$ where $a^0\in
(J(\lambda))^{\aleph_0}$ and
$\varphi_2(x)=(\varphi_0(x),\varphi_1(x),a^0)$ is also embedding.

Since $\varphi_0(F_1)\subseteq [0,\frac{1}{2}]$ and
$\varphi_0(V_4)=1$, then $F_1\cap V_4=\emptyset$. Hence
$(X\setminus (V_3\cup V_4))\cap F_1=F_1\setminus V_3$. Since
$\varphi_0(V_3)=0$, then $V_3\subseteq F_1$. Similarly
$(X\setminus (V_3\cup V_4))\cap F_2=F_2\setminus V_4$ and
$V_4\subseteq F_2$.

\medskip

{\it Construction of $f_1$ and $f_2$} where

$f_1: F_1\rightarrow ((0,\frac{1}{2}]\times \{d_0\}\cup \{0\}\cup
((0,1]\times (S\setminus \{d_0\})))\times
J(\lambda)^{\aleph_0}\times J(\lambda)^{\aleph_0}$,

$f_2: F_2\rightarrow (([\frac{1}{2},1]\times \{d_0\})\cup
\{0\})\times J(\lambda)^{\aleph_0}\times J(\lambda)^{\aleph_0}$
such that $f_1\upharpoonright (F_1\setminus
V_3)=\varphi_2\upharpoonright (F_1\setminus V_3)$ and
$f_2\upharpoonright (F_2\setminus V_4)=\varphi_2\upharpoonright
(F_2\setminus V_4)$.

The mappings $f_1$ and $f_2$ are constructed similarly, but the
image of $f_2$ is more complex, so we will construct only the
mapping $f_2$.

Note that $D_4\subseteq V_4\subseteq F_2$. Let $D_4=\bigcup \{T_i:
i\in \omega\}$ where $T_i\cap T_j=\emptyset$ for $i\neq j$ and
$|T_i|=\lambda$ for each $i,j\in \omega$. The family
$\{F_2\setminus V_4, T_i, i\in\omega\}$ is discrete. Define the
mapping $\varphi: (F_2\setminus V_4)\cup D_4\rightarrow
\mathbb{N}$ by $\varphi((F_2\setminus V_4)\cup T_1)=1$ and
$\varphi(T_i)=i$ for $i>1$. Let $\varphi^*: F_2\rightarrow
\mathbb{R}$ be a continuous extension of $\varphi$.

Let $\Phi_1=(\varphi^*)^{-1}((-\infty, \frac{3}{2}])$ and
$\Phi_i=(\varphi^*)^{-1}([i-\frac{1}{2},i+\frac{1}{2}])$ for
$i>1$. Then the family $\{\Phi_i : i\in\omega\}$ is a closed
locally-finite cover of $F_2$ and $\Phi_i\setminus \bigcup
\{\Phi_j: j\in \omega$, $j\neq i\}\supseteq T_i$ for each $i\in
\omega$.

Let $a=(a_j)\in J(\lambda)^{\aleph_0}=\{b=(b_j): b_j\in
J(\lambda)$ for $j\in \omega\}$. Denote by
$J(\lambda)^{\aleph_0}[a]=\{b=(b_j)\in J(\lambda)^{\aleph_0}:
|j\in \omega : b_j\neq a_j|<\aleph_0\}$.

Choose $a^i=(a^i_j)\in  J(\lambda)^{\aleph_0}$ such that
$J(\lambda)^{\aleph_0}[a^k]\cap
J(\lambda)^{\aleph_0}[a^m]=\emptyset$ for $k\neq m$,
$k,m\in\omega$.

Let $\Phi_0=F_2\setminus V_4$.

{\it Construction of continuous mappings $\psi_k$}

$\psi_k: \bigcup \{\Phi_i: i=0,...,k\}\rightarrow
((\frac{1}{2},1]\times J(\lambda)^{\aleph_0}\times
J(\lambda)^{\aleph_0})\cup \varphi_2(F_1\cap F_2)$ for each
$k=0,1,...$ such that

1). $\psi_k$ is an injection,

2). $\psi_{k+1}$ a extension of $\psi_k$,

3). $((\frac{1}{2},1]\times J(\lambda)^{\aleph_0}\times
J(\lambda)^{\aleph_0}\setminus
\bigcup\{J(\lambda)^{\aleph_0}[a^i]: i=k,...\})\cup
\varphi_2(F_1\cap F_2)\subseteq \psi_{k+1}(\bigcup\{\Phi_i:
i=0,...,k\})\subseteq ((\frac{1}{2},1]\times
J(\lambda)^{\aleph_0}\times J(\lambda)^{\aleph_0}\setminus
\bigcup\{J(\lambda)^{\aleph_0}[a^i]: i=k+1,...\})\cup
\varphi_2(F_1\cap F_2)$,

4). $\psi_{k+1}\upharpoonright (\bigcup\limits_{i=0}^{k+1}
\Phi_i\setminus \bigcup\limits_{i=0}^{k} \Phi_i)$ is a
homeomorphism for each $k\in \omega$.

Let $\psi_0=\varphi_2\upharpoonright \Phi_0$. Assume that
$\psi_0,...,\psi_m$ are constructed.

%The set $\psi_m(\bigcup\limits_{i=0}^m \Phi_i)$ is Borel. Then
%$C_0=\{(\frac{1}{2},1]\times J(\lambda)^{\aleph_0}\times
%(J(\lambda)^{\aleph_0}\setminus
%\bigcup\{J(\lambda)^{\aleph_0}(a^i): i=m+1,...\})\setminus
%\psi_m(\bigcup\limits_{i=0}^m \Phi_i)\}$ is Borel, too.

Let $C_0=\{(\frac{1}{2},1]\times J(\lambda)^{\aleph_0}\times
(J(\lambda)^{\aleph_0}\setminus
\bigcup\{J(\lambda)^{\aleph_0}[a^i]: i=m+1,...\})\setminus
\psi_m(\bigcup\limits_{i=0}^m \Phi_i)\}$.

Since $\lambda=\lambda^{\aleph_0}$ then $|C_0|=\lambda$. Let
$\eta_{m+1}:T_{m+1}\rightarrow C_0$ be a condensation.

For $[\frac{1}{2},1]\times J(\lambda)^{\aleph_0}\times
J(\lambda)^{\aleph_0}$ we write as $I_0\times
\prod\limits_{i=1}^{\infty} J^1_i\times
\prod\limits_{i=1}^{\infty} J^2_i$.

Let $\pi_i^j$ be the projection $[\frac{1}{2},1]\times
J(\lambda)^{\aleph_0}\times J(\lambda)^{\aleph_0}$ onto $J_i^j$
for $j=1,2$ and $i\in \omega$ and $\pi_0$ be the projection
$[\frac{1}{2},1]\times J(\lambda)^{\aleph_0}\times
J(\lambda)^{\aleph_0}$ onto $I_0$.

The set $\bigcup \{\Phi_i: i=0,...,m\}\cup T_{m+1}$ is closed in
$\bigcup \{\Phi_i: i=0,...,m+1\}$. Then $\bigcup \{\Phi_i:
i=0,...,m+1\}\setminus (\bigcup \{\Phi_i: i=0,...,m\}\cup
T_{m+1})=\bigcup\{M_j: j\in\omega\}$ where $M_j\subseteq M_{j+1}$
and $M_j$ is closed in $\bigcup \{\Phi_i: i=0,...,m+1\}$ for each
$j\in \omega$.

Define $f^2_j: (\bigcup \{\Phi_i: i=0,...,m\}\cup T_{m+1}\cup
M_j)\rightarrow J^2_j$ for $j\in \omega$ as follows

$$ f^2_j(x)= \left\{
\begin{array}{lcr}
\pi^2_j\psi_m(x)$, \ \ \ \ if  $x\in \bigcup\limits_{i=0}^m \Phi_i, \\
\pi^2_j\eta_{m+1}(x)$,  \  if  $x\in T_{m+1}, \\
a^j_{m+1}$ , \ \ \ \ \ \ \  if $x\in M_j.\\
\end{array}
\right.
$$

Let $f^{2,*}_j: \bigcup\limits_{i=0}^{m+1} \Phi_i \rightarrow
J^2_j$ be a continuous extension of $f^2_j$.

Fix a homeomorphic embedding $\Omega_j: M_j\rightarrow
\prod\limits_{i=1}^{\infty} J^1_i$ for each $j\in \omega$.

Define $f^1_j:(\bigcup\limits_{i=0}^m \Phi_i\cup T_{m+1}\cup
M_j)\rightarrow J^1_j$ for $j\in \omega$ as follows

$$ f^1_j(x)= \left\{
\begin{array}{lcr}
\pi^1_j\psi_m(x)$, \ \ \ \ if  $x\in \bigcup\limits_{i=0}^m \Phi_i, \\
\pi^1_j\eta_{m+1}(x)$,  \  if  $x\in T_{m+1}, \\
\pi^1_j \Omega_j(x)$ , \ \ \ \ \ \ \  if $x\in M_j.\\
\end{array}
\right.
$$

Let $f^{1,*}_j: \bigcup\limits_{i=0}^{m+1} \Phi_i \rightarrow
J^1_j$ be a continuous extension of $f^1_j$.

Consider $f^0: \bigcup\limits_{i=0}^{m+1} \Phi_i\cup
T_{m+1}\rightarrow I_0$ as follows

$$ f^0(x)= \left\{
\begin{array}{lcr}
\pi_0 \psi_m(x)$, \ \ \ \ if  $x\in \bigcup\limits_{i=0}^m \Phi_i, \\
\pi_0 \eta_{m+1}(x)$,  \  if  $x\in T_{m+1}.\\
\end{array}
\right.
$$

By properties (1-4),  and the definition of $\eta_{m+1}$,
$(f^0)^{-1}(\frac{1}{2})=F_1\cap F_2$.

Let $f^{0,*}: \bigcup\limits_{i=0}^{m+1} \Phi_i \rightarrow I_0$
be a continuous extension of $f^0$ such that
$(f^{0,*})^{-1}(\frac{1}{2})=F_1\cap F_2$.

Define $\psi_{m+1}: \bigcup\limits_{i=0}^{m+1} \Phi_i \rightarrow
 (((\frac{1}{2},1]\times \{d_0\})\times J(\lambda)^{\aleph_0}\times
J(\lambda)^{\aleph_0})\cup \varphi_2(F_1\cap F_2)$ as follows
$\psi_{m+1}(x)=f^{0,*}(x)\times
\{f^{1,*}_j(x)\}_{j=1}^{\infty}\times
\{f^{2,*}_j(x)\}_{j=1}^{\infty}$.

Note that the system of functions $\psi_0, ..., \psi_{m+1}$ is
such that properties (1-4) are satisfied.

Let the function $\psi_k$ is constructed for each $k\in \omega$.

Define $f_2: F_2 \rightarrow (([\frac{1}{2},1]\times \{d_0\})\cup
\{0\})\times J(\lambda)^{\aleph_0}\times J(\lambda)^{\aleph_0}$ as
follows $f_2(x)=\psi_n(x)$ where $n=\min \{k: x\in \Phi_k\}$.

Because $\{\Phi_i\}$ is a closed locally-finite family, $f_2$ is a
continuous function.

By properties (1) and (2), $f_2$ is an injection. By property (3),
$f_2$ is a surjection.

Let $f_1$ and $f_2$ are constructed. Then define $f: X \rightarrow
(J(\lambda)\times J(\lambda)^{\aleph_0}\times
J(\lambda)^{\aleph_0})\cong J(\lambda)^{\aleph_0}$ as follows
$f(x)=f_i(x)$, if $x\in F_i$, $i=1,2$. This definition is correct
because  $F_1\cap F_2\subseteq (F_1\setminus V_3)\cap
(F_2\setminus V_4)$ and, hence, $f_1\upharpoonright (F_1\cap
F_2)=\varphi_2\upharpoonright (F_1\cap F_2)=f_2\upharpoonright
(F_1\cap F_2)$.

Because $f_1$ and $f_2$ are injections, then $f$ is an
injection,too.

Since $f\upharpoonright (X\setminus (V_3\cup
V_4))=\varphi_2\upharpoonright (X\setminus (V_3\cup V_4))$, then
the property (I) is satisfied.

Claim that the property (III) is satisfied.

Since $F_1\setminus V_3\supseteq V_1$ and $f\upharpoonright
(F_1\setminus V_3)=f_1 \upharpoonright (F_1\setminus V_3)=
\varphi_2\upharpoonright (F_1\setminus V_3)$, then
$f(V_1)=\varphi_2(V_1)\subseteq \varphi_0(V_1)\times
J(\lambda)^{\aleph_0}\times J(\lambda)^{\aleph_0}=\{0\}\times
J(\lambda)^{\aleph_0}\times J(\lambda)^{\aleph_0}$. On the other
side, $f(V_4)\subseteq f(F_2)=f_2(F_2)\subseteq
[\frac{1}{2},1]\times J(\lambda)^{\aleph_0}\times
J(\lambda)^{\aleph_0}$. It follows that $\overline{f(V_1)}\cap
\overline{f(V_4)}=\emptyset$.

Similarly, it is proved that $\overline{f(V_3)}\cap
\overline{f(V_2)}=\emptyset$.

Then the condensation $f: X \rightarrow (J(\lambda)\times
J(\lambda)^{\aleph_0}\times J(\lambda)^{\aleph_0})$, satisfying
properties (I) and (III), is constructed.

To construct the function $g$, it suffices to change  $V_1$ with
$V_3$ and $V_2$ with $V_4$ in the construction $f$.

Let $f$ and $g$ are constructed. Define $\tau_1$ ($\tau_2$) as an
initial structure on a set $X$ generated by map $f$ ($g$,
respectively), i.e. $\tau_1$ ($\tau_2$) is the coarsest topology
on $X$ making $f$ ($g$, respectively) continuous.

Note that $(X\setminus (V_3\cup V_4))\cup (X\setminus (V_1\cup
V_2))=X$,  $\tau\upharpoonright (X\setminus (V_3\cup
V_4)=\tau_1\upharpoonright (X\setminus (V_3\cup V_4))$ and
$\tau\upharpoonright (X\setminus (V_1\cup
V_2)=\tau_2\upharpoonright (X\setminus (V_1\cup V_2))$.

In order to prove that $\tau=\sup\{\tau_1,\tau_2\}$, it is
necessary to show that $X\setminus (V_3\cup V_4)$ and $X\setminus
(V_1\cup V_2)$ are closed in $\sup\{\tau_1,\tau_2\}$. It is enough
to prove that $V_i$ ($i=1,2,3,4$) is an open set in
$\sup\{\tau_1,\tau_2\}$.

Since $V_1\subseteq  X\setminus (V_3\cup V_4)$ then, by property
(I), there is an open set $U$ in $J(\lambda)\times
J(\lambda)^{\aleph_0}\times J(\lambda)^{\aleph_0}$ such that
$V_1=(X\setminus (V_3\cup V_4))\cap f^{-1}(U)$. By property (III),
$X\setminus f^{-1}(\overline{f(V_4)})\supseteq V_1$. By property
(IV), $X\setminus g^{-1}(\overline{g(V_3)})\supseteq V_1$. Then
$V_1=(X\setminus (V_3\cup V_4))\cap f^{-1}(U)\supseteq (X\setminus
f^{-1}(\overline{f(V_4)})\cap (X\setminus
g^{-1}(\overline{g(V_3)})\cap f^{-1}(U)\supseteq V_1\cap
f^{-1}(U)$. Hence, $V_1=(X\setminus f^{-1}(\overline{f(V_4)})\cap
(X\setminus g^{-1}(\overline{g(V_3)})\cap f^{-1}(U)$, i.e., $V_1$
is an open set in $\sup\{\tau_1,\tau_2\}$.

Similarly, we can prove that $V_i$ is an open set in
$\sup\{\tau_1,\tau_2\}$ for $i=2,3,4$. Thus
$\tau=\sup\{\tau_1,\tau_2\}$. In particular, the space $(X,\tau)$
admits a condensation onto $J(\lambda)^{\aleph_0}$ where
$J(\lambda)=(J(\lambda), \tau_{\rho})$ is the metrizable hedgehog
of spininess $\lambda$.

Since $\tau_{\rho}$ is the supremum of two topologies $\tau_{d_1}$
and $\tau_{d_2}$, where $(J(\lambda),\tau_{d_i})$ is a Hausdorff
compact space  for $i=1,2$, then $X$ admits a bijective continuous
mapping onto a Hausdorff compact space
$(J_{d_1}(\lambda))^{\aleph_0}$ where
$J_{d_1}(\lambda)=(J(\lambda),\tau_{d_1})$.

\end{proof}

\begin{corollary} {\it Every metric space of weight $\lambda=\lambda^{\aleph_0}$ admits a
bijective continuous mapping onto a Banach space of weight
$\lambda$.}
\end{corollary}

\begin{corollary} {\it Every metric space of weight $\lambda=\lambda^{\aleph_0}$ admits a
bijective continuous mapping onto a Hausdorff compact space.}
\end{corollary}

\begin{corollary} {\it If $X$ is a metric space of
weight $\lambda$ and $D$ is a discrete space of cardinality
$\lambda^{\aleph_0}$ then $X\times D$  admits a condensation onto
a Banach space of weight $\lambda^{\aleph_0}$.}
\end{corollary}

In \cite{bp}, it is proved that every Banach space of weight
$\mathfrak{c}$ admits a condensation onto the Hilbert cube.

Since $\mathfrak{c}=\mathfrak{c}^{\aleph_0}$ then we have the
following results.

\begin{theorem}{\it If $X$ is a metric space of
weight $\mathfrak{c}$ then $X$ admits a condensation onto the
Hilbert cube.}
\end{theorem}

\begin{theorem} {\it If $X$ is a metric space of
weight $\lambda\leq\mathfrak{c}$ and $D$ is a discrete space of
cardinality $\mathfrak{c}$ then $X\times D$  admits a condensation
onto the Hilbert cube.}
\end{theorem}

\bibliographystyle{model1a-num-names}
\bibliography{<your-bib-database>}

\end{document}